\documentclass[10pt,technote]{IEEEtran}

\ifCLASSINFOpdf
\else
\fi
\usepackage[caption=false]{subfig}
\usepackage[dvips]{graphicx}
\usepackage{amsmath, amssymb}
\usepackage{psfrag}
\usepackage{pgfplots}
\usepackage{arydshln}
\usepackage{mathtools}
\usepackage{rotating}

\usepackage{algorithm}
\usepackage{algorithmic}

\DeclareMathOperator*{\minimize}{minimize}
\DeclareMathOperator*{\maximize}{maximize}

\DeclareMathOperator*{\subject}{subject\ to}

\DeclareMathOperator*{\diag}{diag}

\newcounter{thm}

\newcounter{remcount}

\newtheorem{theorem}[thm]{Theorem}

\newtheorem{rem}[remcount]{Remark}
\newtheorem{pf}{Proof}

\usepackage{color}

\begin{document}
%
\title{Distributed, scalable and gossip-free consensus optimization with application to data analysis}
%
%
%
\vspace{-2mm}
\author{Sina~Khoshfetrat~Pakazad, 
        Christian~A.~Naesseth, 
        Fredrik~Lindsten, 
        Anders~Hansson,~\IEEEmembership{Member,~IEEE}%
\thanks{S.\ Khoshfetrat Pakazad is with ..., Sweden. Email: sina.khoshfetrat@gmail.com}
\thanks{C.\ A.\ Naesseth and A.\ Hansson are with the Division of Automatic Control, Department of Electrical Engineering, Link\"oping University, Sweden. Email: \{christian.a.naesseth, anders.g.hansson\}@liu.se.}
\thanks{F.\ Lindsten is with the Division of Systems and Control, Department of Information Technology, Uppsala University, Sweden. Email: fredrik.lindsten@it.uu.se.}}%
\maketitle
\vspace{-2mm}
\begin{abstract}
Distributed algorithms for solving additive or consensus optimization problems commonly rely on first-order or proximal splitting methods. These algorithms generally come with restrictive assumptions 
and at best enjoy a linear convergence rate. Hence, they can require many iterations or communications among agents to converge. In many cases, however, we do not seek a highly accurate solution for consensus problems. Based on this we propose a controlled relaxation of the coupling in the problem which allows us to compute an approximate solution, where
the accuracy of the approximation 
can be controlled by the level of relaxation.
The relaxed problem can be efficiently solved in a distributed way using a combination of primal-dual interior-point methods (PDIPMs) and message-passing.
This algorithm purely relies on second-order methods and thus requires far fewer iterations and communications to converge. This is illustrated in numerical experiments, showing its superior performance compared to existing methods.
\end{abstract}
\vspace{-1mm}
\begin{IEEEkeywords}
Distributed optimization, data analysis, consensus, primal-dual method, data analysis.
\end{IEEEkeywords}

%
\IEEEpeerreviewmaketitle

\vspace{-4mm}
\section{Introduction}
Many optimization problems in e.g.\ machine learning, control and signal processing \cite{boyd:11,has:03,cev:14,ohl:14} can be formulated as 
\begin{align}\label{eq:SumOriginal}
\minimize_{x} \quad \frac{1}{N} \sum_{i = 1}^N f_i(x) + g_i(x),
\end{align}
where the convex functions $f_i: \mathbb R^p \rightarrow \mathbb R$ and $g_i: \mathbb R^p \rightarrow \mathbb R$ are smooth and non-smooth, respectively, such that $F_i = f_i + g_i$ is Lipschitz continuous. Here we assume that $x \in \mathbb R^p$, where $p$ is not overly large, whereas $N$ can be potentially large.

It is sometimes impossible to solve these problems using centralized optimization algorithms. This is commonly due to \emph{computational issues}, e.g., when $N$ is very large, or due to \emph{privacy requirements}. 
In these cases, a solution is provided by distributed algorithms, which solve the optimization problem using a network of computational agents that can collaborate and communicate with one another. These algorithms commonly rely on consensus formulations of the problem and are typically based on first-order splitting methods, see e.g., \cite{ber:97,ned:10,cev:14,shi:15,boyd:11,com:11}. Thus, these methods are slow with sub-linear or at best linear convergence rates, \cite{cev:14,shi:15,ned:10}. Moreover, they sometimes require further assumptions, such as smoothness or strong convexity of the cost function, and are commonly sensitive to the scaling of the problem. Among these algorithms the ones based on proximal point methods or proximal method of multipliers are less sensitive to scaling, see e.g., \cite{eck:89,boyd:11,com:11}, and require less iterations to converge. However, each iteration is generally far more computationally demanding than that of gradient-based or subgradient-based algorithms.

In order to reduce the required number of iterations and sensitivity to scaling, attempts have been made to combine first and second order methods, see e.g., \cite{kho:15b,wei:13,mok:16}. 
However due to reliance on first-order methods, these algorithms still require many iterations to converge. In a distributed setting, this means that they require many communications among computational agents, which gives rise to non-negligible communication overhead. In order to address all these issues, we set out to devise algorithms that solely rely on second-order methods, which allow for \emph{(i)} more efficient handling of the problem data, \emph{(ii)} convergence in few iterations, and \emph{(iii)} efficient use of parallel/distributed computational platforms and cloud/edge computing for solving the problem. To reach this goal, we face two main hurdles. Firstly, using second-order methods for solving \eqref{eq:SumOriginal} in a distributed or parallelized manner is generally not possible due to the fact that the subproblems in \eqref{eq:SumOriginal} are all fully coupled. Secondly, second-order methods cannot be applied directly for solving this problem, due to the fact that the functions $F_i$ are non-differentiable or not continuously differentiable.

In this paper, we show how these hurdles can be overcome. In Section~\ref{sec:Relaxed}, we first present a controlled relaxation to the coupling among the subproblems in~\eqref{eq:SumOriginal}. 
By solving the relaxed formulation we obtain an approximate solution for \eqref{eq:SumOriginal}, where the accuracy of the approximation is controlled by the level of relaxation. This is motivated by the fact that in many applications we do not seek an exact solution to \eqref{eq:SumOriginal}, for instance due to the uncertainty present in the problem data. A similar relaxation was considered in \cite{kop:16}, but based on a different motivation and and for handling streamed data.

Next, we propose to use primal-dual interior-point methods (PDIPMs) (reviewed in Section~\ref{sec:PD}) for solving the relaxed problem. The convergence of these methods is well-established, see e.g., \cite{boyd:04, wri:97}. The proposed relaxation allows us to impose a coupling structure on the problem that can be represented as a tree. This opens up for using message-passing to compute the search directions and other parameters \emph{exactly} through conducting recursions over the tree structure (see also \cite{kho:16}). Hence, distributing the computations does not jeopardize the convergence of the PDIPM. Message-passing is briefly discussed in Section~\ref{sec:DC}. The resulting algorithm provides superior performance in comparison to recently developed high-performance distributed algorithms, as shown in Section~\ref{sec:NR} using multiple numerical experiments.
\section{A Controlled Relaxation and Reformulation}\label{sec:Relaxed}
In order to impose a desirable structure on \eqref{eq:SumOriginal}, let us consider a relaxation of the problem given as

\vspace{-4mm}
\small
\begin{subequations}\label{eq:Relaxed}
\begin{align}
\minimize_{x, x^i}& \quad \frac{1}{N}\sum_{i=1}^{N}  F_i(x^i) \label{eq:Relaxed-a}\\
\subject & \quad \| x - x^i \|^2 \leq \varepsilon^2, \quad i = 1, \dots, N. \label{eq:Relaxed-c}
\end{align}
\end{subequations}
\normalsize
Note that the terms in \eqref{eq:Relaxed-a} are decoupled. The approximate coupling among these terms are now described using the constraints in \eqref{eq:Relaxed-c}. Let us denote an optimal solution of the problem in \eqref{eq:SumOriginal} with $x^\ast$ and that of~\eqref{eq:Relaxed} with $x_{\textrm{rel}}^\ast$ and $x^{i,\ast}$. It is possible to compute satisfactory suboptimal solutions for \eqref{eq:SumOriginal} by solving \eqref{eq:Relaxed}, as quantified by the following theorem.
\begin{theorem}\label{thm:thm1}
Let us assume that the Lipschitz constant for each $F_i$ is denoted by $L_i$. Then we have

\vspace{-4mm}
\small
\begin{align}\label{eq:Suboptimality}
\frac{1}{N}\sum_{i = 1}^{N} F_i(x_{\textrm{rel}}^*) - F_i(x^{\ast}) &\leq \frac{\varepsilon}{N} L,
\end{align}
\normalsize
where $L = \sum_{i = 1}^{N} L_i$. Furthermore, if the cost function is strongly convex with modulus $m$, we have
\(
\| x_{\textrm{rel}}^* - x^*\|^2 \leq \frac{2\varepsilon L}{Nm}.
\)
\end{theorem}
\begin{pf}
See Appendix \ref{app:appB}. \hfill$\square$
\end{pf}

If the tolerated suboptimality of the solution is $\varepsilon_{\textrm{tol}}$, choosing $\varepsilon = N\varepsilon_{\textrm{tol}}/L$ guarantees that $x_{\textrm{rel}^*}$ gives a satisfactory solution. Moreover, if the problem is strongly convex, given a threshold $\varepsilon_{\textrm{var}}$ concerning the accuracy of the solution, if we choose $\varepsilon = \frac{Nm\varepsilon_{\textrm{var}}}{2L}$, we can guarantee that the obtained solution will satisfy the accuracy requirements. It goes without saying that the smaller the $\varepsilon$, the more accurate the computed solution. However, care must be taken as choosing extremely small values for this parameter can give rise to numerical issues. In general, e.g., provided that the data is normalized, we can compute accurate enough solutions using moderately small values of $\varepsilon$, see Section \ref{sec:NR}.


\begin{algorithm}[tb]
	\caption{\small Primal-dual Interior-point Method \cite{wri:97,boyd:04}}\label{alg:PDIPM}
	\small
	\begin{algorithmic}[1]
		\STATE{Given feasible iterates}
		\REPEAT
		\STATE{Compute the primal-dual search directions}
		\STATE{Compute appropriate primal and dual step sizes}
		\STATE{Update primal and dual iterates}
		\STATE{Update the perturbation parameter}
		\UNTIL{stopping criteria is satisfied}
	\end{algorithmic}
\end{algorithm}

In this paper, we devise a distributed algorithm for solving the relaxed problem \eqref{eq:Relaxed}, purely relying on second-order methods.
Due to non-smoothness of the objective function, however, second-order algorithms cannot be directly applied. Instead, we introduce additional variables and constraints in order to equivalently reformulate the problem as (see e.g., \cite{boyd:04}),%
%

\vspace{-4mm}
\begin{subequations}\label{eq:SumIP}
\small
\begin{align}
\minimize_{x,x^i,t^i}& \quad \frac{1}{N}\sum_{i = 1}^{N} h_i(x^i,t^i)\\
\subject & \quad G^i(x^i,t^i) \leq 0, \quad i = 1, \dots, N\label{eq:SumIP-b}\\
& \quad A^i\begin{bmatrix} x^i\\t^i \end{bmatrix} = b^i, \quad i = 1, \dots, N\label{eq:SumIP-c} \\
&\quad \| x - x^i \|^2 \leq \varepsilon^2, \quad i = 1, \dots, N\label{eq:SumIP-d}
\end{align}
\normalsize
\end{subequations}
where $h_i: \mathbb R^{p+d_i} \rightarrow \mathbb R$ is smooth, the variables $t^i \in \mathbb R^{d_i}$ denote the additional variables, and the constraints in~\eqref{eq:SumIP-b} and~\eqref{eq:SumIP-c} are the additional inequality and equality constraints, with $A^i \in \mathbb R^{u_i \times (p+d_i)}$ and $G^i : \mathbb R^{p+d_i} \rightarrow \mathbb R^{m_i}$. For instance when the non-smooth terms in the objective function are indicator functions for convex sets, the problem can be reformulated by removing these terms and adding the corresponding constraints to the problem. Another common approach for reformulating e.g., problems of the form in~\eqref{eq:Relaxed} as in \eqref{eq:SumIP}, is through the use of epigraph reformulations \cite{boyd:04}.
Other approaches are also possible, see e.g. \cite{boyd:04}, but for the sake of brevity we do not discuss them here.


%
%
%
%
%

\section{Primal-dual Interior-point Methods}\label{sec:PD}
PDIPMs are the state-of-the-art iterative solvers for problems like \eqref{eq:SumIP}.
A generic description of these methods is given in Algorithm \ref{alg:PDIPM}, \cite{wri:97,boyd:04}.
Let us denote the dual variables for the constraints in \eqref{eq:SumIP-b}, \eqref{eq:SumIP-c} and \eqref{eq:SumIP-d} with $z^i$, $v^i$ and $\lambda_i$, respectively. At each iteration $k$, given feasible primal and dual iterates, i.e., such that $G^i(x^{i,(k)},t^{i,(k)}) < 0$, $\| x^{(k)} - x^{i,(k)} \|^2 < \varepsilon^2$, $z^{i,(k)} > 0$ and $\lambda^{(k)}_i > 0$, one way of computing the search directions requires solving an equality constrained quadratic program, see~\cite{kho:16}. Particularly, for the problem in \eqref{eq:SumIP}, the QP that needs to be solved takes the form%

\vspace{-4mm}
\small
\begin{subequations}\label{eq:SumIPQP}
\begin{align} \notag
\minimize_{\Delta x,\Delta x^i,\Delta t^i} &\quad \sum_{i = 1}^{N}\frac{1}{2}
\left[
\begin{array}{@{}c@{}} \Delta x^i \\ \Delta t^i \\\hdashline \Delta x \end{array} \right]^T
\underbrace{\left[ \begin{array}{@{}c:c@{}} H^{i,(k)}_{ll} & H^{i,(k)}_{lg}\\ \hdashline (H^{i,(k)}_{lg})^T & H^{i,(k)}_{gg} \end{array}\right]}_{H_{\textrm{pd}}^{i,(k)}}
\left[\begin{array}{@{}c@{}} \Delta x^i \\ \Delta t^i \\\hdashline \Delta x  \end{array}\right] \\
 &\hspace{3cm}+ \left[\begin{array}{@{}c@{}} \Delta x^i \\ \Delta t^i \\\hdashline \Delta x  \end{array}\right]^T \left[\begin{array}{@{}c@{}} r^{i,(k)}_{l} \\ \hdashline r^{i,(k)}_{g}  \end{array}\right]\\
\subject & \quad A^i \begin{bmatrix} \Delta x^i \\ \Delta t^i\end{bmatrix} = r_{\textrm{primal}}^{i,(k)}
\end{align}
\end{subequations}
\normalsize
Through solving \eqref{eq:SumIPQP} we can compute the primal directions $\Delta x^{i,(k+1)}$, $\Delta t^{i,(k+1)}$ and $\Delta x^{(k+1)}$ together with the dual directions $\Delta v^{i,(k+1)}$, see \cite{kho:16}. It is then possible to compute the remaining dual variables' directions
$\Delta z^{i,(l+1)}$ and $\Delta \lambda_i^{(k+1)}$; the explicit expressions are provided in Appendix~\ref{app:app1} together with expressions for the data matrices appearing in \eqref{eq:SumIPQP}.
For more details on how these matrices are formed, see e.g., \cite{wri:97} or \cite[Sec. 5 and 6]{kho:16}.

At this point, we can compute an appropriate step size, e.g., using back-tracking line search that assures feasibility of the iterates and persistent decrease of the norm of primal and dual residuals, \cite{wri:97,noe:06,boyd:04}. We can then update the iterates and the procedure is continued until certain stopping criteria are satisfied. These are commonly based on primal and dual residuals norms and the so-called surrogate duality gap, see \cite{wri:97,kho:16} for more details.
%
%

During the run of a PDIPM, the main computational burden arises from the computation of the search directions, which requires solving \eqref{eq:SumIPQP}. Indeed, the cost of this can be prohibitive in many cases. Also, for problems that come with privacy requirements, the computations cannot be done in a centralized manner. However, due to the coupling structure of the problem in \eqref{eq:SumIP}, which is also inherited by \eqref{eq:SumIPQP}, it is possible to distribute the computations at each iteration of the PDIPM using message-passing (or dynamic programming) over trees as discussed in \cite{kho:16}.
 Next, we show how this can be done for the problem under study.
\section{Distributed Computations}\label{sec:DC}
\begin{figure}[t]
	\begin{center}
		\includegraphics[width=3cm]{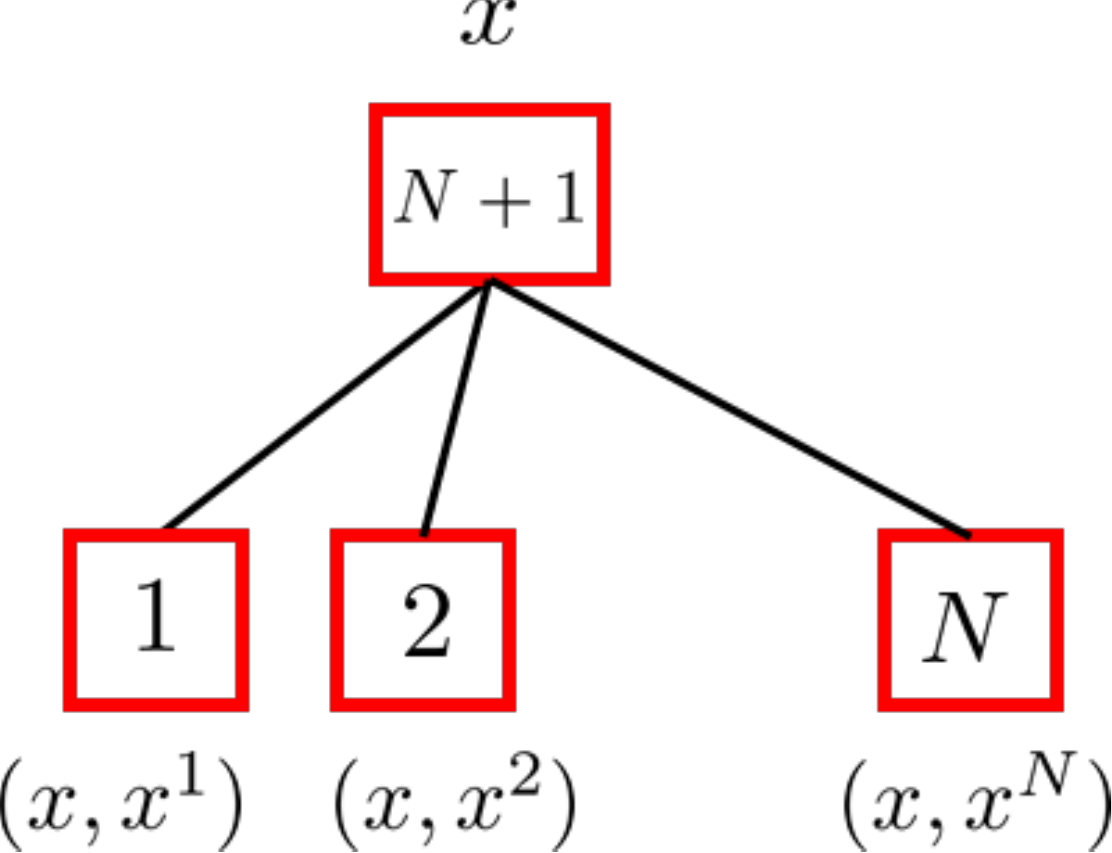}    
		\caption{\small Tree representation of the coupling structure of \eqref{eq:Relaxed}.}
		\label{fig:Large}
		\vspace{-5mm}
	\end{center}
\end{figure}

Let us reconsider the problem in \eqref{eq:SumIP}. This problem is made up of $N$ subproblems, each of which is defined by a term in the cost function and its corresponding constraint set described by each term in \eqref{eq:SumIP-b}--\eqref{eq:SumIP-d}. The coupling structure of this problem can be represented using the tree illustrated in Figure \ref{fig:Large}.

Recall that the main computational burden
 of a PDIPM applied to \eqref{eq:SumIP}
corresponds to the computations of the search directions, i.e.\ solving the QP in \eqref{eq:SumIPQP}. Note that this QP inherits the same tree representation and, hence, we can solve the problem in \eqref{eq:SumIPQP} by conducting message-passing upward and downward through the tree, see \cite{kho:16,kol:09}. For this purpose, we first assign each subproblem $i$ in \eqref{eq:SumIPQP} to each leaf $i$ of the tree. Then, considering the star-shape structure of the tree, at each iteration $k$, agents at the leaves of the tree compute their messages to the root of the tree simultaneously. Specifically each agent first computes the search directions for the local variables $\Delta x^i$ and $\Delta v^i$ as a function of $\Delta x$ by solving

\vspace{-4mm}
\small
\begin{align*}
\left[ \begin{array}{@{}c:c@{}} H^{i,(k)}_{ll} & (A^i)^T\\ \hdashline A^i & 0 \end{array}\right] \left[\begin{array}{@{}c@{}} \Delta x^i \\ \Delta t^i \\\hdashline \Delta v^i  \end{array}\right] = \left[ \begin{array}{@{}c@{}} -r^{i,(k)}_{l} \\ \hdashline r_{\textrm{primal}}^{i,(k)}   \end{array}\right] - \left[ \begin{array}{@{}c@{}} H^{i,(k)}_{lg} \\ \hdashline 0  \end{array}  \right]\Delta x,
\end{align*}
\normalsize
the result of which can be written compactly as

\vspace{-4mm}
\small
\begin{align}\label{eq:localIter}
\left[\begin{array}{@{}c@{}} \Delta x^i \\ \Delta t^i \\\hdashline \Delta v^i  \end{array}\right] = \left[ \begin{array}{@{}c@{}} u_1^{i,(k)} \\ \hdashline u_2^{i,(k)}   \end{array}\right] + \left[ \begin{array}{@{}c@{}} U^{i,(k)}_{1} \\ \hdashline U^{i,(k)}_{2}  \end{array}  \right]\Delta x.
\end{align}
\normalsize
By inserting the solutions from \eqref{eq:localIter} into the cost function of the local subproblems, we obtain quadratic functions in $\Delta x$, with Hessian and linear terms given by

\vspace{-4mm}
\small
\begin{align*}
Q^{i,(k)} & = H^{i,(k)}_{gg} + (U^{i,(k)}_{1})^TH^{i,(k)}_{ll}U^{i,(k)}_{1} +\notag\\& \quad \quad \quad \quad \quad \quad (U^{i,(k)}_{1})^TH^{i,(k)}_{lg} + (H^{i,(k)}_{lg})^TU^{i,(k)}_{1},\\
q^{i,(k)} & = r^{i,(k)}_{g} + (U^{i,(k)}_{1})^Tr^{i,(k)}_{l} + (H^{i,(k)}_{lg})^T u_1^{i,(k)} +\notag\\& \quad \quad \quad \quad \quad \quad \quad \quad \quad \quad \quad  \quad (U^{i,(k)}_{1})^TH^{i,(k)}_{ll}u_1^{i,(k)},
\end{align*}
\normalsize
respectively. These quadratic functions are then sent to the root.
The agent at the root 
will then form and solve the optimization problem
\begin{align*}
\minimize_{\Delta x} \quad \frac{1}{2} \sum_{i = 1}^N \Delta x^T Q^{i,(k)} \Delta x + \Delta x^T q^{i,(k)}
\end{align*}
which gives the search direction $\Delta x$, that it then communicates downwards to its children. Each agent at the leaves of the tree, having received $\Delta x$, can then compute its local variables' search directions using \eqref{eq:localIter}.

Notice that computing stepsizes and residuals, updating the perturbation parameters and checking the termination condition require conducting summing or computing minimum or maximum of local quantities over the tree, and hence, we can use the same computational structure for this purpose, see \cite[Sec. 6.3]{kho:16}. Consequently, combining message-passing and PDIPMs results in a scalable and distributed algorithm, that purely relies on second-order methods, for solving the problem in \eqref{eq:Relaxed}. A generic summary of the proposed algorithm is given in Algorithm \ref{alg:DPDA}. Note that mixing message-passing and PDIPMs does not affect their convergence and the proposed method thus inherits properties such as superlinear convergence and finite termination of PDIPMs. See Appendix \ref{app:appC} for further discussion.
%
 \begin{algorithm}[tb]
\caption{\small Distributed Primal-dual Algorithm (DPDA)}\label{alg:DPDA}
\small
\begin{algorithmic}[1]
\small
\STATE{Given $k = 0$, $\mu>1$, $\varepsilon>0$, $\epsilon_{\text{d}}>0$, $\epsilon_{\text{feas}}>0$, $x^{(k)}$, $x^{i,(0)}, t^{i,(0)}$ $z^{i,(0)}$, $v^{i,(0)}$ and $\lambda^{(0)}_i$, such that $G^i(x^{i,(0)},t^{i,(0)}) < 0$, $\| x^{(0)} - x^{i,(0)} \|^2 < \varepsilon^2$, $z^{i,(0)}, \lambda^{(0)}_i > 0$ for all $i = 1, \dots, N$,~$\hat \eta^{(0)}$ and $\delta = \mu m /\hat \eta^{(0)}$}
\REPEAT
\STATE{Perform message-passing upwards and downwards through the tree in Figure \ref{fig:Large} to compute the search directions}
\STATE{Compute a proper step size, $\alpha^{(k+1)}$, by performing upward-downward passes through the tree, see \cite[Sec. 6.3]{kho:16} for details.}
\STATE{Update the primal and dual iterates using the computed search directions and step size}
\STATE{Perform upward-downward pass through the tree to decide whether to terminate the algorithm and/or to update the perturbation parameter $\delta = \mu m /\hat \eta^{(k+1)}$.}
\STATE {$k = k + 1$.}
\UNTIL{the algorithm is terminated}
\normalsize
\end{algorithmic}
\normalsize
\end{algorithm}
\begin{rem}
The coupling structure presented in Figure \ref{fig:Large} is imposed based on the way the constraints in~\eqref{eq:Relaxed-c} are introduced. It is possible to impose other structures, e.g., chain-like or balanced trees, by modifying the way these constraints are introduced. Doing so requires recomputing the bounds calculated in Section \ref{sec:Relaxed} to match this structure.
\end{rem}

Note that all agents have access to their local variables updates $x^{i,(k)}$ and that of the global one $x^{(k)}$ as depicted in Figure \ref{fig:Large}, and all agents consider the solution for $x$ as the computed parameters. This means that we have exact consensus among agents.
\vspace{-2mm}
\section{Numerical Experiments}\label{sec:NR}
In this section we apply the proposed algorithm DPDA to robust least squares and logistic regression problems, and compare its performance with that of alternating direction method of multipliers (ADMM), \cite{boyd:11}, and algorithms presented in \cite{shi:15} and \cite{mok:16}. We refer to these algorithms as EXTRA and ESOM, respectively. These algorithms are chosen based on their superior performance in comparison to commonly used algorithms for distributedly solving problems of the form~\eqref{eq:SumOriginal}. We compare the performance of the algorithms based on their iterations count and computational time. We do not claim that any of the algorithms (including DPDA) has been implemented in their most efficient manner, which can potentially affect the reported computational time, whereas the iteration count to be less susceptible to this. Another reason for considering the iteration count is that it corresponds to the number of communications among agents. This is a good performance measure since for many existing algorithms the communication overhead is the most significant bottleneck which can create significant latency, see e.g., \cite{jag:14,li:14}.
\vspace{-2mm}
\subsection{Robust Least Squares Problem}
We apply DPDA to a least squares problem given as
\begin{align}\label{eq:RLS}
\minimize_{x} \quad \sum_{i = 1}^N \sum_{j = 1}^{n_i} \phi_M(A_j^i x - Y_j^i),
\end{align}
where $A^i \in \mathbb R^{n_i\times p}$ with $A^i_j$ denoting the $j$th row of $A^i$, and $\phi: \mathbb R \rightarrow \mathbb R$ is the Huber penalty function defined~as
\begin{align*}
\phi_M(u) = \begin{cases} u^2 \quad\quad\quad\quad\quad\ \  |u|\leq M \\ M(2|u| - M) \ \ \ |u| > M \end{cases}.
\end{align*}
\normalsize
We assume that each agent $i$ has access to its own measurements $Y^i = A^i x + e^i$ with $A^i \in \mathbb R^{n_i \times p}$ and $e^i \sim \mathcal N(0, \sigma^2I)$ is the measurement noise. In this experiment, $N = 10$, $n_i = 20$ and $p = 10$, the matrices $A^i$ have been generated randomly based on a uniform distribution in the interval $[0, 1]$ and the parameters $x$ used for producing the data have been also generated randomly in the interval $[0, 20]$.

\begin{figure}[tb]
	\begin{center}
		\includegraphics[width=4.2cm]{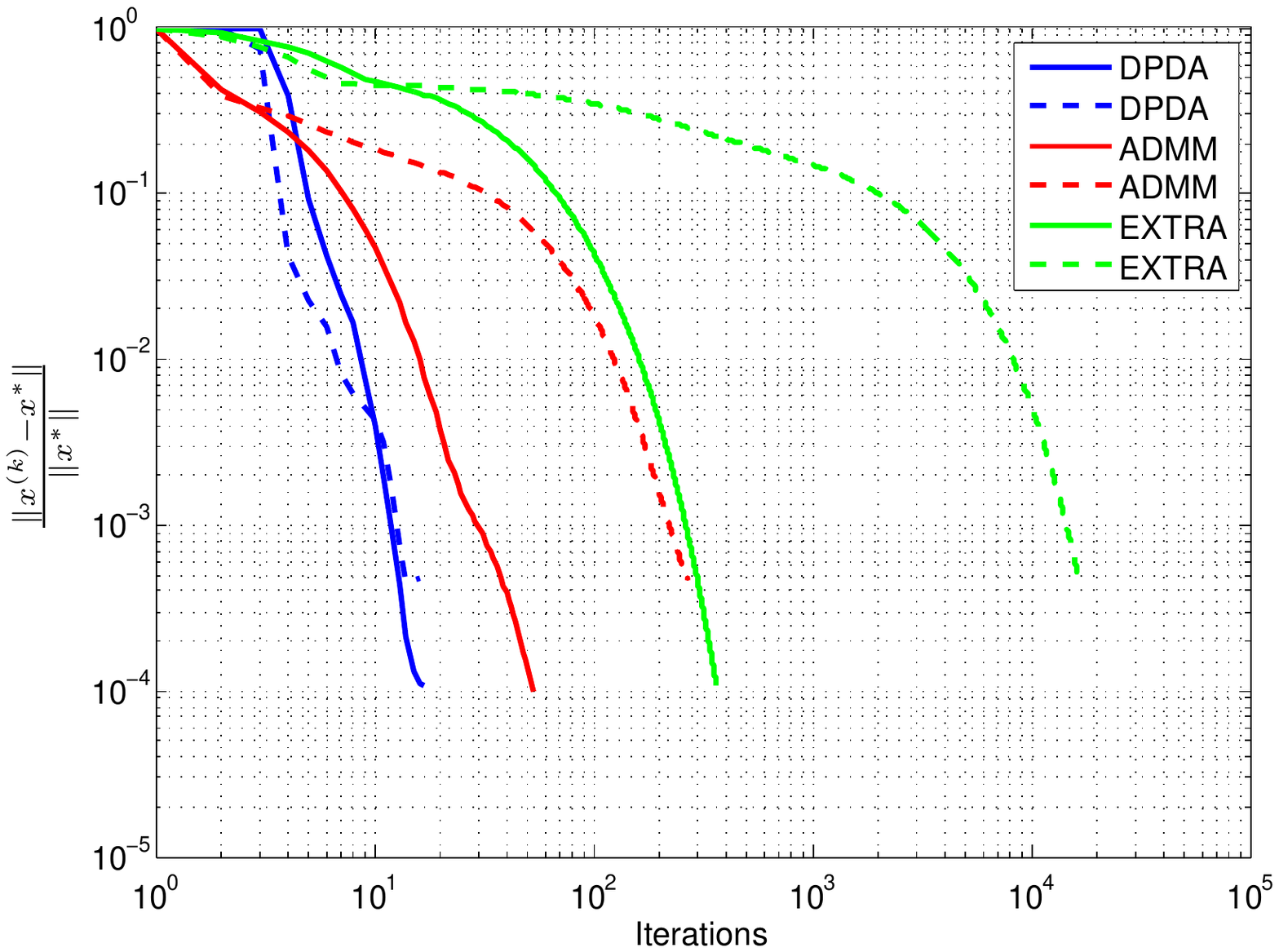}
		\includegraphics[width=4.2cm]{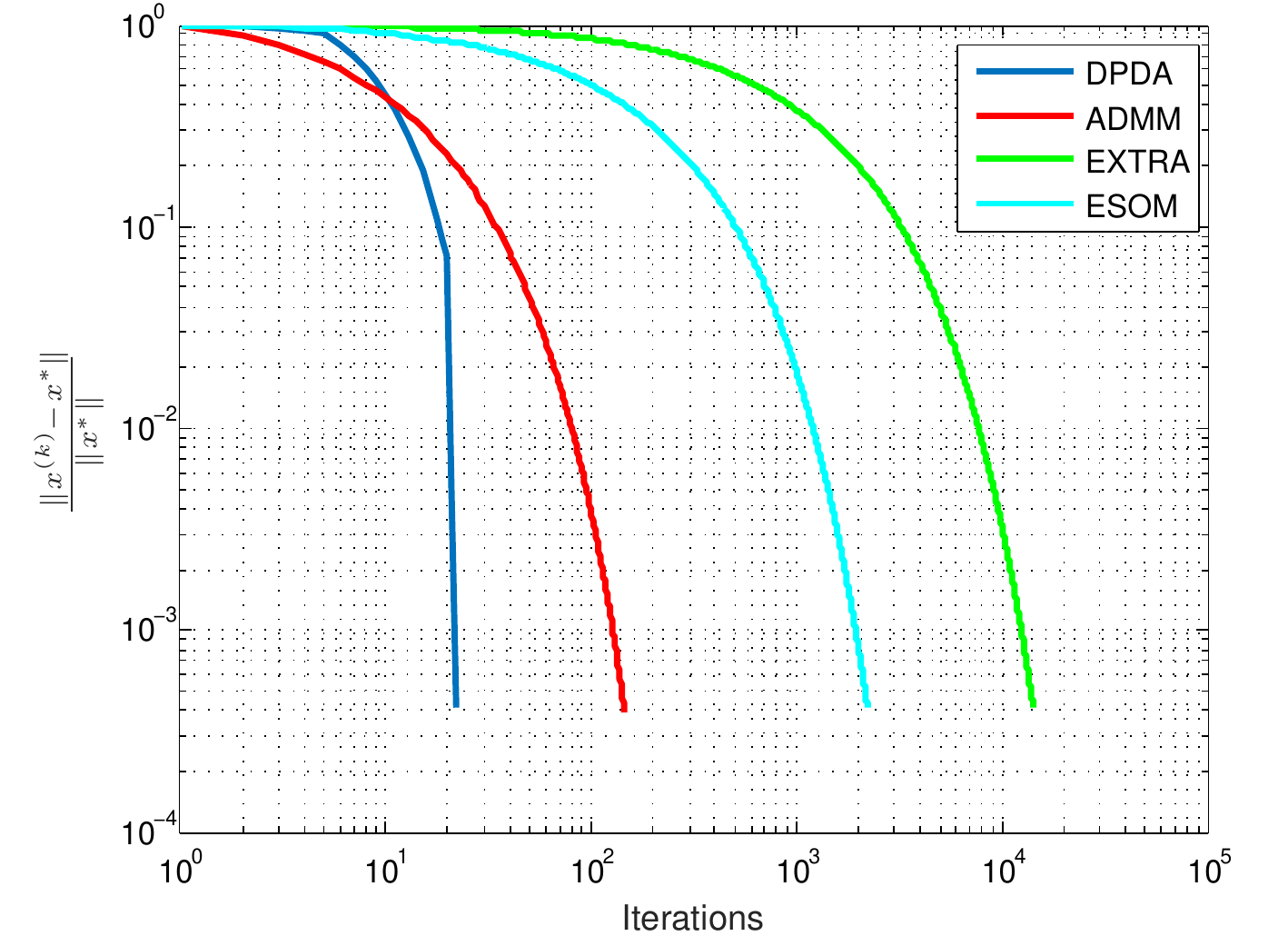}    
		\vspace{-2mm}
		\caption{\small \emph{(Left)} Results based on a robust least squares problem where the condition number of $A$ is $6.56$, depicted using solid lines and where the condition number of $A$ is $56.92$ depicted using the dashed lines.
			\emph{(Right)} Results based on a logistic regression problem.}
		\label{fig:RobustLS1}
		\vspace{-5mm}
	\end{center}
\end{figure}

Notice that although the cost function for this problem is smooth, it is not twice continuously differentiable. This means ESOM cannot be applied to this problem and in order to use DPDA we
use the equivalent reformulation (see \cite{boyd:04})

\vspace{-4mm}
\small
\begin{equation*}
\begin{split}
\minimize_{x,u^i,v^i} &\ \ \sum_{i = 1}^N \|u^i\|^2 + M\mathbf 1^T v^i\\
\subject & \ \begin{rcases*} -u^i-v^i \leq A^ix - Y^i \leq u^i + v^i \\\quad \quad \quad \  0 \leq u^i \leq M\mathbf 1 \\\quad \quad \quad  v^i \geq 0\end{rcases*}, i = 1, \dots, N.
\end{split}
\end{equation*}
\normalsize
 Note that ESOM cannot be used for solving this formulation either, as it includes constraints. We set $M = 1$ and the relaxation level for DPDA to $\varepsilon = 10^{-3}$. The algorithm parameters are chosen as $\mu = 10$, $\beta = 0.4$ and $\alpha = 0.1$. The matrices $W$ and $\tilde W$ used in EXTRA are chosen in the same way as in \cite[Sec. 4]{shi:15} for the star-shaped graph in Figure~\ref{fig:Large}. The other parameters in this algorithm, and those in ADMM, are tuned manually to maximize performance.

 The results are shown in Figure \ref{fig:RobustLS1} (left). The figure reports results for two experiments. First for a coefficient matrix $A$ generated as described above (condition number $6.56$).
 Second, in order to study the effect of the scaling of the problem, for a coefficient matrix $A$ obtained by manipulating its singular values to increase the condition number by almost a factor ten (to $56.92$).
The optimal solution for~\eqref{eq:RLS}, $x^*$, has been computed using CVX, \cite{gra:08}. From these figures, we observe that DPDA requires far fewer iterations than ADMM and EXTRA, and hence communications among agents, to converge to a solution. Furthermore, it is much less sensitive to the scaling of the problem. Despite this, one should note that since EXTRA is a gradient-based method, the computational complexity of each of its iterations is linear in the number of variables. This is in contrast to DPDA where the computational cost of each of iterations is cubic in the number of variables. Moreover, ADMM has a higher computational complexity EXTRA and DPDA. This is because at each iteration of the ADMM, each agent needs to solve a robust least squares problem for updating its local variables. Due to this, ADMM in fact has the worst per-iteration computational complexity among these algorithms. The computational time for this experiment are reported in Table \ref{tab:RLS} (the first two rows). As can be seen from the table, EXTRA and DPDA provide comparable performance, and DPDA clearly outperforms EXTRA for worse conditioned problems. Having said that, DPDA will potentially perform worse on problems with large number of features. However, as was mentioned in the introduction these problems are not the focus of this paper.
    \begin{table}
  \centering
   \caption{\small Results for a robust least squares problem, first row for cond($A$) = 6.56 and second row for cond($A$) = 56.92. Results for a logistic regression problem in third row.}
   \small
   \vspace{-2mm}
\begin{tabular}{|c|c|c|c|c|}
   \hline
   &DPDA & ADMM & ESOM & EXTRA \\ \cline{1-5}
   Time [sec]& 6.19 &  365.17& -- & 5.65 \\ \cline{1-5}
   Time [sec]& 5.85 &  649.96& -- & 146.85 \\ \cline{1-5}
       \hline
       \hline
          Time [sec] & 6.92 & 1342.01 & 1340.26 & 84.45 \\ \cline{1-5}
       \hline
\end{tabular}
\label{tab:RLS}
\end{table}

\vspace{-2mm}
\subsection{Logistic Regression Problem}

We continue our investigation of the performance of DPDA by conducting an experiment based on a logistic regression problem. A logistic regression problem can be written as

\vspace{-4mm}
\small
\begin{align}\label{eq:LR}
\maximize_{x} \quad \sum_{i = 1}^N \left[\sum_{j = 1}^{n_i} \left(Y^i_j\Phi^i_jx - \log(1 + e^{\Phi^i_jx})\right) + \frac{\rho}{N}\| x \|^2\right]
\end{align}
\normalsize
where $\Phi^i \in \mathbb R^{n_i \times p}$ with $\Phi^i_j$ as its $j$th row, and $Y^i_j \in \{0, 1\}$. The regularization term $\frac{\rho}{N}\| x \|^2$ is generally added to prevent over-fitting to the data, where $\rho>0$ is the so-called penalty or regularization parameter. This parameter has been chosen as $\rho = 1$. The data for this problem concerns the classification problem of radar returns from the ionosphere and has been taken from \cite{lic:13}. For this problem $p = 34$ and we have considered $350$ data points, that we assume are divided among $N = 10$ agents. This means that $n_i = 35$. The results from these experiments are illustrated in Figure \ref{fig:RobustLS1} (right).
The computational times are reported in Table \ref{tab:RLS} (the third row).
Similar to the previous experiment we see that DPDA clearly outperforms the other algorithms, and ADMM has the highest per-iteration computational complexity. This is because each agent at each iteration needs to solve an optimization problem similar to \eqref{eq:LR} in order to update its local variables. 

\section{Conclusions}
In this paper we proposed a distributed PDIPM for computing approximate solutions for convex consensus problems. This was done by first proposing a relaxed formulation of the consensus problem. Solving this problem results in an approximate solution of the consensus problem, where we showed how the accuracy of the computed solution can be controlled by the relaxation level. The imposed coupling structure in the relaxed problem enabled us to distribute the computations of each iteration of a PDIPM using message-passing. We showed the performance of the proposed algorithm using numerical experiments based on robust least squares and logistic regression problems, using both synthetic and real data. In this paper we did not discuss a standard approach for choosing $\varepsilon$. We plan to address this as a future line of research, possibly through introduction of efficient methodologies for global or local scaling of the problem data.

\ifCLASSOPTIONcaptionsoff
  \newpage
\fi



%
\bibliographystyle{plain}
\bibliography{IEEETrans}
\appendices
\section{Description of the data matrices in \eqref{eq:SumIPQP}} \label{app:app1}

Here we present a detailed description of the data matrices defining \eqref{eq:SumIPQP}. Let us start with the Hessian matrix of the cost function that is given as

\small
\begin{multline*}
H^{i,(k)}_{\textrm{pd}} = \left[\begin{array}{c:c} \nabla^2 h_i^{(k)} & 0\\ \hdashline 0 & 0 \end{array}\right] + \sum_{j = 1}^{m_i} z^{i,(k)}_j \left[\begin{array}{c:c} \nabla^2 G_j^{i,(k)} & 0\\ \hdashline 0 & 0  \end{array}\right] - \\ \sum_{j = 1}^{m_i} \frac{z^{i,(k)}_j}{G^{i,(k)}_j} \left[\begin{array}{c} \nabla_{x^i}G^{i,(k)}_j \\ \nabla_{t^i}G^{i,(k)}_j \\ \hdashline 0 \end{array}\right]\left[\begin{array}{c} \nabla_{x^i}G^{i,(k)}_j \\ \nabla_{t^i}G^{i,(k)}_j \\ \hdashline 0 \end{array}\right]^T + \\2\lambda^{(k)}_i \left[\begin{array}{cc:c} I & 0 & -I\\0 & 0 & 0\\ \hdashline -I & 0 & I \end{array}\right] - \frac{2\lambda_i^{(k)}}{\|x^{(k)} -x^{i,(k)}\|^2-\varepsilon^2}\times\\\left[\begin{array}{c} x^{i,(k)}-x^{(k)}\\ 0\\ \hdashline x^{(k)} - x^{i,(k)} \end{array}\right]\left[\begin{array}{c} x^{i,(k)}-x^{(k)}\\ 0\\ \hdashline x^{(k)} - x^{i,(k)} \end{array}\right]^T,
\end{multline*}
\normalsize
with $\nabla^2 h_i^{(k)} = \begin{bmatrix} \nabla_{x^ix^i}h^{(k)}_i & \nabla_{x^it^i}h^{(k)}_i \\  \star & \nabla_{t^it^i}h^{(k)}_i\end{bmatrix}$ and the matrices $\nabla^2 G_j^{i,(k)}$ are defined similarly. The coefficient vector defining the linear term in the cost function can be extracted as below

\small
\begin{multline*}
\left[\begin{array}{c} r^{i,(k)}_{l} \\ \hdashline r^{i,(k)}_{g}  \end{array}\right] = r_{\textrm{dual}}^{i,(k)} + \left[\begin{array}{c} (D_{x^i}G^{i,(k)})^T \\ (D_{t^i}G^{i,(k)})^T \\ \hdashline 0 \end{array} \right]\diag(G^{i,(k)})^{-1} r^{i,(k)}_{\textrm{cent}} + \\\frac{2r^{i,(k)}_{\text{Q}}}{\| x^{(k)} - x^{i,(k)} \|^2 - \varepsilon^2}\left[\begin{array}{c} x^{i,(k)}-x^{(k)} \\0 \\ \hdashline x^{(k)}-x^{i,(k)} \end{array} \right],
\end{multline*}
\normalsize
with $D_{\star}G^{i,(k)} = \begin{bmatrix} \nabla_{\star} G^{i,(k)}_1 & \dots & \nabla_{\star} G^{i,(k)}_{m_i} \end{bmatrix}^T$ and

\small
\begin{align*}
r_{\textrm{dual}}^{i,(k)} &= \left[\begin{array}{c} \nabla_{x^i}h^{(k)}_i \\ \nabla_{t^i}h^{(k)}_i \\ \hdashline 0 \end{array} \right] + \sum_{j = 1}^{m_i}z_j^{i,(k)} \left[\begin{array}{c} \nabla_{x^i}G_j^{i,(k)} \\ \nabla_{t^i}G_j^{i,(k)} \\ \hdashline 0 \end{array} \right] + \\ & \quad\quad\quad 2 \lambda^{(k)}_i \left[\begin{array}{c} x^{i,(k)}-x^{(k)} \\0 \\ \hdashline x^{(k)}-x^{i,(k)} \end{array} \right] + \left[\begin{array}{c} (A^i)^T \\ \hdashline 0 \end{array} \right]v^{i,(k)},\\
r^{i,(k)}_{\textrm{cent}} &= -\diag(z^{i,(k)})G^{i,(k)} - \frac{1}{\delta}\mathbf 1,\\
r_{\textrm{Q}}^{i,(k)} &= -\frac{\lambda^{(k)}_i}{\|x^{(k)}-x^{i,(k)}\|^2-\varepsilon^2} - \frac{1}{\delta},
\end{align*}
\normalsize
where for the sake of notational ease we have denoted a function $G(x)$ evaluated at $x^{(k)}$ with $G^{(k)}$. Here $\delta$ is referred to as the perturbation parameter. Having computed the primal variables' directions $\Delta x^{i,(k+1)}$, $\Delta t^{i,(k+1)}$ and $\Delta x^{(k+1)}$ together with the dual variables' directions $\Delta v^{i,(k+1)}$ by solving~\eqref{eq:SumIPQP}, we can compute the remaining dual variables' directions as

\small
\begin{subequations}\label{eq:SumIPQP2}
	\begin{align}\notag
	\Delta z^{i,(l+1)} &= -\diag( G^{i,(k)})^{-1} \left( \diag(z^{i,(k)}) \times \right. \\&   \left.\begin{bmatrix}D_{x^i} G^{i,(k)} &D_{t^i} G^{i,(k)} \end{bmatrix} \begin{bmatrix} \Delta x^{i,(k+1)} \\ \Delta t^{i,(k+1)} \end{bmatrix} - r^{i,(k)}_{\text{cent}} \right)\\\notag
	\Delta \lambda_i^{(k+1)} &= \frac{1}{\|x^{(k)} -x^{i,(k)}\|^2-\varepsilon^2}\times \\&   \left( \lambda_i^{(k)}\begin{bmatrix} x^{i,(k)}-x^{(k)}\\ x^{(k)} - x^{i,(k)} \end{bmatrix}^T \begin{bmatrix} \Delta x^{i,(k+1)} \\ \Delta x^{(k+1)} \end{bmatrix}- r^{i,(k)}_{\textrm{Q}} \right)
	\end{align}
\end{subequations}
\normalsize
Given $\mu > 1$ and once we have updated the primal and dual variables, the perturbation parameter can then be updated as $\delta = \mu m/\hat \eta^{(k+1)}$ with $m = N + \sum_{i = 1}^N m_i$ and
\begin{multline*}
\hat \eta^{(k+1)} = \sum_{i=1}^N -\lambda_i^{(k+1)}\left( \| x^{i,(k+1)}-x^{(k+1)} \|^2 - \varepsilon^2 \right)-\\(z^{i,(k+1)})^TG^{i,(k+1)}
\end{multline*}
\normalsize
denoting the surrogate duality gap.
\section{Proof of Theorem \ref{thm:thm1}}\label{app:appB}
Firstly, notice that we have
\small
\begin{align*}
\underbrace{\frac{1}{N}\sum_{i = 1}^{N}  F_i(x^{i,\ast})}_{L_b} \leq \frac{1}{N}\sum_{i = 1}^{N} F_i(x^*) \leq \underbrace{\frac{1}{N}\sum_{i = 1}^{N}  F_i(x_{\textrm{rel}}^*)}_{U_b},
\end{align*}
\normalsize
where the first inequality follows from the fact that \eqref{eq:Relaxed} is a relaxation of \eqref{eq:SumOriginal} and the second inequality follows from the fact that $x^*$ is optimal for \eqref{eq:SumOriginal} but $x_{\textrm{rel}}^*$ is not. Then we have

\small
\begin{multline*}
\sum_{i = 1}^{N} F_i(x_{\textrm{rel}}^*) - F_i(x^{\ast}) \leq N (U_b-L_b)
\leq \sum_{i = 1}^{N} \| F_i(x_{\textrm{rel}}^*) - F_i(x^{i,\ast}) \|\\
\leq \sum_{i = 1}^{N} L_i \| x_{\textrm{rel}}^* - x^{i,\ast} \|
\leq \varepsilon L.
\end{multline*}
\normalsize
which proves \eqref{eq:Suboptimality}. Under the assumption that the cost function is strongly convex, we have
\small
\begin{align*}
\| x_{\textrm{rel}}^* - x^*\|^2 \leq \frac{2}{Nm}\sum_{i = 1}^{N} F_i(x_{\textrm{rel}}^*) - F_i(x^*)\leq \frac{2\varepsilon L}{Nm}
\end{align*}
\normalsize
which completes the proof.
\section{Convergence properties of the proposed method}\label{app:appC}
PDIPMs have been shown to converge and that they enjoy favorable convergence properties such as superlinear convergence and finite termination, see \cite[Ch. 6 and 7]{wri:97} for a full discussion and technical presentation. Algorithm \ref{alg:DPDA} is obtained by distributing the computations of each iteration of the PDIPM. As was shown in \cite[Sec. 6.1 and Thm 6.4]{kho:16}, the computed search directions using message-passing are exact. Furthermore the exact computation of the parameters can be trivially distributed within a message-passing framework, see \cite[Sec. 6.3]{kho:16}. This means that the distributed computations do not jeopardize the convergence of the PDIPM and does not affect its convergence properties.
%

%
%
%




\end{document}